\documentclass{amsart}

\usepackage{amsmath,amssymb,amscd,amsfonts}

\newtheorem{theorem}{Theorem}
\newcommand{\bt}{\begin{theorem}}
\newcommand{\et}{\end{theorem}}
\newtheorem{lemma}{Lemma}
\newcommand{\bl}{\begin{lemma}}
\newcommand{\el}{\end{lemma}}
\newcommand{\beq}{\begin{equation}}
\newcommand{\eeq}{\end{equation}}
\newcommand{\benum}{\begin{enumerate}}
\newcommand{\eenum}{\end{enumerate}}
\newcommand{\N}{\ensuremath{ \mathbf N }}
\newcommand{\mcb}{\ensuremath{ \mathcal B}}
\newcommand{\mcs}{\ensuremath{ \mathcal S}}

\begin{document}

\title[sequences without geometric progressions]
{Irrational numbers associated to sequences without geometric progressions}
\author{Melvyn B. Nathanson}
\address{Department of Mathematics\\
Lehman College (CUNY)\\
Bronx, NY 10468}
\email{melvyn.nathanson@lehman.cuny.edu}

\author{Kevin O'Bryant}
\address{Department of Mathematics\\
College of Staten Island (CUNY)\\
Staten Island, NY 10314}
\email{kevin@member.ams.org}

\subjclass[2010]{11B05 11B25, 11B75, 11B83, 05D10, 11J99.}
\keywords{Geometric progression-free sequences, Ramsey theory, irrationality.}

\date{\today}

\begin{abstract}
Let $s$ and $k$ be integers with $s \geq 2$ and $k \geq 2$ .  
Let $g_k^{(s)}(n)$ denote the cardinality of the largest subset of the set 
$\{1,2,\ldots, n\}$ that contains no geometric progression of length $k$ 
whose common ratio is a power of $s$.   
Let $r_k(\ell)$ denote the cardinality of the largest subset of the set 
$\{0,1,2,\ldots, \ell -1\}$ that contains no arithmetric progression of length $k$ .  
The limit 
\[
\lim_{n\rightarrow \infty} \frac{g_k^{(s)}(n)}{n}
=  (s-1)  \sum_{m=1}^{\infty} 
 \left( \frac{1}{s} \right)^{   \min \left( r_k^{-1}(m)\right) } 
\]
exists and converges to an irrational number.  
\end{abstract}

\maketitle

\section{Maximal subsets without geometric progressions}
Let \N\ denote the set of positive integers.  
For every real number $x$, the \emph{integer part} of $x$, denoted $[x]$, 
is the unique integer $n$ such that $n \leq x < n+1$.

Let $s \geq 2$ be an integer.  
Every positive integer $a$ can be written uniquely in the form 
\[
a = bs^v
\]
where $b$ is a positive integer not divisible by $s$ and $v$ is a nonnegative integer.  
If $G$ is a finite geometric progression of length $k$ 
whose common ratio is a power of $s$, say, $s^d$, 
then 
\[
G = \{a\left(s^d\right)^j: j = 0,1,\ldots, k-1  \}.
\]
Writing $a$ in the form $a = bs^v$, we have
\beq    \label{NoGP:form}
G = \{bs^{v+dj}: j = 0,1,\ldots, k-1  \} \subseteq \{bs^i: i \in \N_0 \}
\eeq
and so the set of exponents of $s$  in the finite geometric progression $G$ 
is the finite arithmetic progression $\{v+dj: j = 0,1,\ldots, k-1  \}$.
Conversely, if $P$ is a finite arithmetic progression 
of $k$ nonnegative integers and if $b$ is a positive integer 
not divisible by $s$, then $\{bs^i :i \in P\}$ is a geometric 
progression of length $k$.  

Let $\ell$ and $k$ be positive integers with $k \geq 2$.  
Let $r_k(\ell)$ denote the cardinality of the largest subset of the set 
$\{0,1,2,\ldots, \ell -1\}$ that contains no arithmetic progression of length $k$.  
Note that $r_k(\ell) = \ell$ for $\ell=1,\ldots, k-1$, 
that $r_k(k) = k-1$, and that, 
for every $\ell \in \N$, there exists $\varepsilon_{\ell} \in \{0,1\}$ such that  
\[
r_k(\ell+1) = r_k(\ell) + \varepsilon_{\ell}.
\]
Thus, the function $r_k:\N \rightarrow \N$ is nondecreasing and surjective.
This implies that, for every positive integer $m$, the set
\[
r_k^{-1}(m) = \{\ell \in \N: r_k(\ell) = m\}
\]
is a nonempty set of consecutive integers, and so
\beq    \label{NoBP:MaxMin}
\max\left(r_k^{-1}(m) \right) + 1 = \min \left(r_k^{-1}(m+1) \right) 
\eeq
and
\beq    \label{NoBP:Max}
\min \left(r_k^{-1}(m) \right) \geq m
\eeq
for all $m \in \N$.  

\bl         \label{NoGP:lemma:szem}
Let  $u_m = \min \left( r_k^{-1}(m) \right)$ for $m \in \N$.
The sequence $(u_m)_{m=1}^{\infty}$ is a strictly increasing sequence 
of positive integers such that 
\[
\limsup_{m\rightarrow \infty} (u_{m+1}-u_m) = \infty.
\]
\el

\begin{proof}
Identity~\eqref{NoBP:MaxMin} implies that 
the sequence $(u_m)_{m=1}^{\infty}$ is strictly increasing.
We use Szemer\' edi's theorem, which states that $r_k(\ell) = o(\ell)$,
to prove that the sequence $(u_m)_{m=1}^{\infty}$ has unbounded gaps.

Note that $u_1 = 1$.  
If $\limsup_{m\rightarrow \infty} (u_{m+1}-u_m) < \infty$, 
then there is an integer $c \geq 2$
such that  $u_{m+1}-u_m < c$ for all $m \in \N$.  
It follows that 
\begin{align*}
\max\left( r_k^{-1}(m)  \right) + 1
& = \min\left( r_k^{-1}(m+1)  \right) \\
& = u_{m+1} \\
& = \sum_{i=1}^m (u_{i+1} - u_{i})+u_1 \\
& < cm+1.
\end{align*}
Thus, $\max\left( r_k^{-1}(m)  \right) < cm$ and so $r_k(cm) > m$.
Equivalently, 
\[
\frac{r_k(cm)}{cm} > c > 0
\]
and 
\[
\liminf_{\ell \rightarrow \infty} \frac{r_k(\ell) }{\ell} \geq c >0.
\]
This contradicts Szemer\' edi's theorem, and completes the proof.
\end{proof}

For $k \geq 2$, let $g_k(n)$ denote the cardinality of the largest subset of the set 
$\{1,2,\ldots, n\}$ that contains no geometric progression of length $k$.  
Rankin~\cite{rank60} introduced this function, and it has been 
investigated by 
M.~Beiglb{\"o}ck, V.~Bergelson, N.~Hindman, and 
D.~Strauss~\cite{beig-berg-hind-stra06}, 
by Brown and Gordon~\cite{brow-gord96}, 
and by Riddell~\cite{ridd69}.  
The best  upper bound for the function $g_k(n)$ is due to 
Nathanson and O'Bryant~\cite{nath-obry13a}.

For $s \geq 2$ and $k \geq 2$, let $g_k^{(s)}(n)$ denote the cardinality 
of the largest subset of the set 
$\{1,2,\ldots, n\}$ that contains no geometric progression of length $k$ 
whose common ratio is a power of $s$.   
We shall prove that the limit 
\beq       \label{NoGP:LimitSeries}
\lim_{n\rightarrow \infty} \frac{g_k^{(s)}(n)}{n}
=  (s-1)  \sum_{m=1}^{\infty} 
 \left( \frac{1}{s} \right)^{   \min \left( r_k^{-1}(m)\right) } 
\eeq
exists and converges to an irrational number.

\section{Maximal geometric progression free sets}

\bl                    \label{NoGP:lemma:formula}
If $k$ and $s$ are integers with $k \geq 2$ and $s \geq 2$, then 
\[
g_k^{(s)}(n) =  \sum_{b\in \mcb_n} r_k\left(1 +  [\log_s(n/b)]  \right).
\]
\el

\begin{proof}
Let $n$ be a positive integer, and let
\[
\mcb_n = \{b\in \{1,2,\ldots, n\} : \text{$s$ does not divide $b$} \}.
\]
If $b \in \mcb_n$ and $i \in \N_0$, then $bs^i \leq n$ if and only if $0 \leq i \leq \log_s(n/b)$.  
We define  
\begin{align*}
T(b) & = \{t \in \{1,2,\ldots, n\} : t = bs^i \text{ for some } i \in \N_0\} \\
& = \{ bs^i: i = 0,1,\ldots, [\log_s(n/b)] \}.
\end{align*}
Then  $b \in T(b)$, and 
\[
\{1,2,\ldots, n\} = \bigcup_{b\in \mcb_n} T(b)
\]
is a partition of $\{1,2,\ldots, n\} $ into pairwise disjoint nonempty subsets.  

If the set $\{1,2,\ldots, n \}$ contains a finite geometric progression of length $k$ 
whose common ratio is a power of $s$, 
then, by~\eqref{NoGP:form}, this geometric progression is a subset of $T(b)$ 
for some $b \in \mcb_n$, and the set of exponents of $s$ is a finite arithmetic 
progression of length $k$ contained in the
set of consecutive integers $\{0,1,\ldots, [\log_s(n/b)] \}$.  
It follows that the largest cardinality of a subset of $T(b)$ 
that contains no $k$-term geometric progression 
is equal to the largest cardinality of a subset of $\{0,1,\ldots, [\log_s(n/b)] \}$ 
that contains no $k$-term arithmetic progression.  This number is 
\[
r_k\left(1 +  [\log_s(n/b)]  \right).
\]
If $A_n$ is a  subset  of $\{1,2,\ldots, n\}$ of maximum cardinality 
that contains no $k$-term geometric progression 
whose common ratio is a power of $s$, then 
\[
|A_n \cap T(b)| = r_k\left(1 +  [\log_s(n/b)]  \right).
\]
Because $A = \bigcup_{b\in \mcb_n} T(b)$ is a partition of $\{1,\ldots, n\}$, 
it follows that 
\[
|A_n| = \sum_{b\in \mcb_n} |A_n \cap T(b)|  = \sum_{b\in \mcb_n} r_k\left(1 +  [\log_s(n/b)]  \right).
\]
This completes the proof. 
\end{proof}

\section{Construction of an irrational number}

\bl        \label{NoGP:lemma}
Let $s$ be an integer with $s \geq 2$.
Let $x$ and $y$ be real numbers with $x < y$.  
The number of integers $n$ such that $x < n \leq y$ 
and $s$ does not divide $n$ is
\[
\left( \frac{s-1}{s} \right) (y-x) + O(1).
\]
\el

\begin{proof}
For every real number $x$, the interval $(x,x+s]$ contains exactly $s$ integers.  
These integers are consecutive, so $(x,x+s]$ contains exactly $s-1$ integers 
not divisible by $s$.  Let $x$ and $y$ be real numbers with $x < y$, 
and let 
\[
h = \left[ \frac{y-x}{s} \right].  
\]
Then
\[
x + hs \leq y < x + (h+1)s.
\]
The interval $(x,x+hs]$ contains exactly $(s-1)h$ integers not divisible by $s$, 
and the  interval $(x,x+(h+1)s]$ contains exactly $(s-1)(h+1)$ integers 
not divisible by $s$.  
If $N$ denote the number of integers in the interval $(x,y]$ 
that are not divisible by $s$, then 
\[
(s-1)h \leq N \leq (s-1)(h+1)
\]
and so 
\[
 \frac{y-x}{s} - 1 < h \leq \frac{N}{s-1} \leq h+1\leq \frac{y-x}{s} +1.
\]
Equivalently, 
\[
\left( \frac{s-1}{s}\right)(y-x) - (s-1) < N \leq \left( \frac{s-1}{s}\right)(y-x) + s-1.
\]
This completes the proof.
\end{proof}

\bt                       \label{NoGP:theorem:series}
Let $k$ and $s$ be integers with $k \geq 2$ and $s \geq 2$.  The limit  
\[
\lim_{n\rightarrow \infty} \frac{g_k^{(s)}(n)}{n}
=  (s-1)  \sum_{m=1}^{\infty} 
 \left( \frac{1}{s} \right)^{   \min \left( r_k^{-1}(m)\right) } 
\]
exists and converges to an irrational number.  
\et

\begin{proof}
For every positive integer $b$ we have 
\[
1 +  [\log_s(n/b)]  = \ell
\]
if and only if 
\[
\frac{n}{s^{\ell}} < b \leq \frac{sn}{s^{\ell}}.
\]
By Lemma~\ref{NoGP:lemma}, the number of  integers in this interval 
that are also in $\mcb_n$, that is, are not divisible by $s$, is 
\[
\left(\frac{s-1}{s} \right)   \frac{(s-1)n}{s^{\ell}} +  O(1)
=   \frac{ n(s-1)^2}{s^{\ell+ 1} }  + O(1).
\]
Because $1 \in \mcb_n$, we have 
\[
L = L(n) = \max\left\{ 1 +  [\log_s(n/b)] : b \in \mcb_n \right\} = 1 + [\log_s n] .
\]
Also, if $\ell \leq L$, then $r_k(\ell) \leq \ell \leq L$.
By Lemma~\ref{NoGP:lemma:formula},
\begin{align*}
g_k^{(s)}(n) 
& = \sum_{b\in \mcb_n} r_k\left(1 +  [\log_s(n/b)]  \right) \\
& = \sum_{\ell = 1}^L r_k\left(\ell  \right) \times |\{b\in \mcb_n: \ell 
= 1 +  [\log_s(n/b)] \}  |\\
& = \sum_{\ell = 1}^L r_k\left(\ell  \right) \left(  \frac{ n (s-1)^2 }{s^{\ell + 1}} 
+ O(1) \right) \\
& = \frac{ n (s-1)^2}{s}  \sum_{\ell = 1}^L \frac{ r_k(\ell ) }{s^{\ell }} 
+  O\left( \sum_{\ell = 1}^L r_k(\ell ) \right) \\
& = \frac{ n (s-1)^2}{s}  \sum_{\ell = 1}^L \frac{ r_k(\ell ) }{s^{\ell }} 
+  O\left(L^2 \right) \\
& = n \left( \frac{(s-1)^2}{s}  \sum_{\ell = 1}^L \frac{ r_k(\ell ) }{s^{\ell }}
+  O\left(  \frac{\log_s^2n}{n} \right)\right) \\
& = n \left(  \frac{(s-1)^2}{s}  \sum_{\ell = 1}^L \frac{ r_k(\ell ) }{s^{\ell }}
+ o(1) \right).
\end{align*}
Let $M = M(n) = r_k(L(n))$.  We have 
\begin{align*}
\sum_{\ell = 1}^L \frac{ r_k(\ell ) }{s^{\ell }} 
& = \sum_{m=1}^{M-1} m\sum_{\ell \in r_k^{-1}(m)}  \frac{ 1 }{s^{\ell}} 
+ M \sum_{\ell \in r_k^{-1}(m)\cap \{1,\ldots, L\} }  \frac{ 1 }{s^{\ell}} \\
& = \sum_{m=1}^{M-1} m\sum_{\ell = \min\left( r_k^{-1}(m) \right)}^{\max \left( r_k^{-1}(m) \right)}   \frac{ 1 }{s^{\ell}} 
+ M \sum_{\ell =  \min\left( r_k^{-1}(m) \right)}^L  \frac{ 1 }{s^{\ell}} \\
& =  \frac{s}{s-1}  \sum_{m=1}^{M-1} m \left( \left( \frac{1}{s} \right)^{  \min\left( r_k^{-1}(m) \right) } -  
 \left( \frac{1}{s} \right)^{  \max \left( r_k^{-1}(m) \right) +1 } \right) \\
& \qquad + \frac{s}{s-1} M \left( \left( \frac{1}{s} \right)^{  \min\left( r_k^{-1}(M) \right) } -  
 \left( \frac{1}{s} \right)^L \right) \\
& =  \frac{s}{s-1}  \sum_{m=1}^{M-1} m \left( \left( \frac{1}{s} \right)^{  \min\left( r_k^{-1}(m) \right) } -  
 \left( \frac{1}{s} \right)^{  \min \left( r_k^{-1}(m+1) \right)} \right) \\
 & \qquad + \frac{s}{s-1} M \left( \left( \frac{1}{s} \right)^{  \min\left( r_k^{-1}(M) \right) } -  
 \left( \frac{1}{s} \right)^L \right) \\
 & =  \frac{s}{s-1} \left(  \sum_{m=1}^M 
m \left( \frac{1}{s} \right)^{   \min \left( r_k^{-1}(m)\right) } 
- \sum_{m=2}^M (m-1) 
\left( \frac{1}{s} \right)^{   \min\left( r_k^{-1}(m)\right) }
- M  \left( \frac{1}{s} \right)^L   \right)  \\
& =  \frac{s}{s-1}  \sum_{m=1}^M 
 \left( \frac{1}{s} \right)^{   \min \left( r_k^{-1}(m)\right) } 
- \frac{sM}{s-1} 
\left( \frac{1}{s} \right)^L \\
& =  \frac{s}{s-1}  \sum_{m=1}^M 
 \left( \frac{1}{s} \right)^{   \min \left( r_k^{-1}(m)\right) } +o(1)
\end{align*}
because 
\[
\frac{sM}{s-1} \left( \frac{1}{s} \right)^L
 \ll \frac{M}{s^L}  \ll \frac{M}{ s^{   \min\left(   r_k^{-1}(M) \right)    } } 
  \ll \frac{M}{s^M} 
\]
by inequality~\eqref{NoBP:Max}.
Therefore,
\begin{align*}
\frac{g_k^{(s)}(n)}{n} 
& =  \frac{(s-1)^2}{s} 
\left(\frac{s}{s-1}  \sum_{m=1}^M 
 \left( \frac{1}{s} \right)^{   \min \left( r_k^{-1}(m)\right) } +o(1)\right) + o(1) \\
 & = (s-1) \sum_{m=1}^M 
 \left( \frac{1}{s} \right)^{   \min \left( r_k^{-1}(m)\right) } +o(1)
\end{align*}
and so 
\[
\lim_{n\rightarrow \infty} \frac{g_k^{(s)}(n)}{n}
=  (s-1) \sum_{m=1}^{\infty} 
 \left( \frac{1}{s} \right)^{   \min \left( r_k^{-1}(m)\right) }.
\]
The infinite series converges to a real number $\theta \in (0,1)$, 
and the ``decimal digits to base $s$'' of $\theta$ are 0 or 1.  
The number $\theta$ is rational if and only if these digits 
are eventually periodic, but Lemma~\ref{NoGP:lemma:szem} 
implies that there are unbounded gaps between successive digits equal to 1.
Therefore, $\theta$ is irrational.  
This completes the proof.  
\end{proof}

\section{Open problems}
\benum
\item
Let $k$ and $s$ be integers with $k \geq 2$ and $s \geq 2$.  
Is  the number
\[
\lim_{n\rightarrow \infty} \frac{g_k^{(s)}(n)}{n}
=  (s-1)  \sum_{m=1}^{\infty} 
 \left( \frac{1}{s} \right)^{   \min \left( r_k^{-1}(m)\right) } 
\]
transcendental?  

\item
Let  $u_m = \min \left( r_k^{-1}(m) \right)$ for $m \in \N$.
Prove that the sequence $(u_m)_{m=1}^{\infty}$ 
is not eventually periodic without using Szemer\' edi's theorem.

\item
Let $s$ and $s'$ be integers with $2 \leq s < s'$.  
Is it true that $g_k^{(s')}(n) \leq g_k^{(s)}(n)$ for all $n \in \N$
and that $g_k^{(s')}(n) < g_k^{(s)}(n)$ for all sufficiently large $n \in \N$?

\item
Let $\mcs$ be a finite set of integers such that $s \geq 2$  for all $s \in \mcs$.  
For $k \geq 2$, let $g_k^{(\mcs)}(n)$ denote the cardinality 
of the largest subset of the set 
$\{1,2,\ldots, n\}$ that contains no geometric progression of length $k$ 
whose common ratio is a power of $s$ for some $s \in \mcs$.   
Does 
\[
\lim_{n\rightarrow \infty} \frac{g_k^{(\mcs)}(n)}{n}
\]
exist?  If so, can this limit be expressed by an infinite series 
analogous to ~\eqref{NoGP:LimitSeries}?

\eenum

\def\cprime{$'$} \def\cprime{$'$} \def\cprime{$'$}
\providecommand{\bysame}{\leavevmode\hbox to3em{\hrulefill}\thinspace}
\providecommand{\MR}{\relax\ifhmode\unskip\space\fi MR }
\providecommand{\MRhref}[2]{%
  \href{http://www.ams.org/mathscinet-getitem?mr=#1}{#2}
}
\providecommand{\href}[2]{#2}


\begin{thebibliography}{1}

\bibitem{beig-berg-hind-stra06}
M.~Beiglb{\"o}ck, V.~Bergelson, N.~Hindman, and D.~Strauss,
  \emph{Multiplicative structures in additively large sets}, J. Combin. Theory
  Ser. A \textbf{113} (2006), no.~7, 1219--1242.

\bibitem{brow-gord96}
B.~E. Brown and D.~M. Gordon, \emph{On sequences without geometric
  progressions}, Math. Comp. \textbf{65} (1996), no.~216, 1749--1754.


\bibitem{nath-obry13a}
M.~B.~Nathanson and K.~O'Bryant, 
\emph{On sequences without geometric progressions} (2013),
available at http://arXiv:1306.0280.

\bibitem{rank60}
R.~A. Rankin, \emph{Sets of integers containing not more than a given number of
  terms in arithmetical progression}, Proc. Roy. Soc. Edinburgh Sect. A
  \textbf{65} (1960/1961), 332--344 (1960/61).

\bibitem{ridd69}
J.~Riddell, \emph{Sets of integers containing no {$n$} terms in geometric
  progression}, Glasgow Math. J. \textbf{10} (1969), 137--146.

\end{thebibliography}
\end{document}